\documentclass[11pt,leqno]{amsart}
\usepackage{graphicx} 
\usepackage{amscd}
\usepackage{color}
\usepackage[symbol]{footmisc}
\usepackage{amssymb}
\usepackage{amsfonts}
\usepackage{latexsym}
\usepackage{verbatim}
\usepackage{bbm}
\usepackage[shortlabels]{enumitem}
\usepackage{amsxtra}
\usepackage{amsopn}
\usepackage{booktabs}
\usepackage{multirow}
\usepackage{array}
\usepackage{color}
\usepackage{amsmath,amsthm,amssymb}
\usepackage{amscd}
\usepackage{amsfonts}
\usepackage{latexsym}
\usepackage{verbatim}
\usepackage{pb-diagram}
\usepackage{color}
\usepackage{tikz}
\usepackage[all]{xy}

\usepackage[backref=page]{hyperref}
\renewcommand*{\backref}[1]{}
\renewcommand*{\backrefalt}[4]{%
	\ifcase #1 (Not cited).%
	\or        (Cited on page~#2).%
	\else      (Cited on pages~#2).%
	\fi}
\hypersetup{colorlinks=true,linkcolor=blue,citecolor=blue,urlcolor=black}
\DeclareMathOperator{\Imm}{Im}
\DeclareMathOperator{\Ker}{Ker}
\newcommand{\del}{\partial}
\newcommand{\delbar}{\overline{\del}}

\newcommand{\C}{\mathbb{C}}

\renewcommand{\epsilon}{\varepsilon}
\renewcommand{\phi}{\varphi}

\theoremstyle{plain}
\newtheorem{thm}{Theorem}
\newtheorem{prop}[thm]{Proposition}
\newtheorem{lemma}[thm]{Lemma}
\newtheorem{cor}[thm]{Corollary}

\theoremstyle{definition}

\newtheorem{rmk}[thm]{Remark}
\newtheorem{es}[thm]{Example}

\title[$\del\delbar$-Lemma and Bott-Chern cohomology of twistor spaces]{$\del\delbar$-Lemma and Bott-Chern cohomology of twistor spaces}

\author{Anna Fino}
\address[Anna Fino]{Dipartimento di Matematica ``G. Peano'', Universit\`{a} degli studi di Torino \\
Via Carlo Alberto 10\\
10123 Torino, Italy\\
\& Department of Mathematics and Statistics, Florida International University\\
Miami, FL 33199, United States}
\email{annamaria.fino@unito.it, afino@fiu.edu}

\author{Gueo Grantcharov}
\address[Gueo Grantcharov]{Department of Mathematics and Statistics \\
Florida International University\\
Miami, FL 33199, United States\\
\& Institute of Mathematics and Iformatics\\
Bulgarian Academy of Sciences\\
8 Acad. G. Bonchev str. 1113, Sofia, Bulgaria}
\email{grantchg@fiu.edu}

\author{Nicoletta Tardini}
\address[Nicoletta Tardini]{Dipartimento di Scienze Matematiche, Fisiche e Informatiche\\
Unit\`{a} di Mate\-matica e Informatica,
Universit\`{a} degli Studi di Parma\\
Parco Area delle Scienze 53/A, 43124 \\
Parma, Italy}
\email{nicoletta.tardini@unipr.it}
\author{Adriano Tomassini}
\address[Adriano  Tomassini]{Dipartimento di Scienze Matematiche, Fisiche e Informatiche\\
Unit\`{a} di Mate\-matica e Informatica,
Universit\`{a} degli Studi di Parma\\
Parco Area delle Scienze 53/A, 43124 \\
Parma, Italy}
\email{adriano.tomassini@unipr.it}

\author{Luigi Vezzoni}
\address[Luigi Vezzoni]{Dipartimento di Matematica ``G. Peano'', Universit\`{a} degli studi di Torino \\
Via Carlo Alberto 10\\
10123 Torino, Italy}
\email{luigi.vezzoni@unito.it}


\begin{document}

\maketitle

\begin{abstract}  In the paper we  study the Bott-Chern and Aeppli cohomologies of the twistor space
of a compact self-dual 4-manifold and we characterize the validity of the $\partial \overline \partial$-lemma. We also   compute explicitly the Dolbeault cohomology of  the twistor space $Z$ of the flat $4$-dimensional torus, which is known to not satisfy the $\partial\overline{\partial}$ lemma.

\end{abstract}

\section{Introduction}

The Bott-Chern and Aeppli cohomology groups of a compact complex manifold $X$:


\[
H^{\bullet,\bullet}_{BC}(X) := \frac{\ker \partial \cap \ker \overline\partial}{\mathrm{im}\, \partial\overline\partial}\;, \qquad
H^{\bullet,\bullet}_{A}(X) := \frac{\ker \partial\overline\partial}{\mathrm{im}\, \partial + \mathrm{im}\, \overline\partial}
\]

provide an important geometric and analytic information when $X$ does not admit any K\"ahler metric.



For instance, in \cite{angella-tomassini} it is shown that a compact complex manifold satisfies the $\partial \overline\partial$-lemma if and only if its Bott-Chern cohomology satisfies a Fr\"olicher-type equality. However, these spaces are very hard to compute and there are only few examples in which any information about them is available.  In this paper we consider the case when $X$ is a twistor space of a compact self-dual 4-manifold.

{\it Twistor spaces} were introduced by R. Penrose  \cite{Penrose} in the 60s as a tool to transform solutions of some nonlinear differential equations in mathematical physics into holomorphic objects in complex geometry. The theory was formalized and generalized in \cite{atiyah-hitchin-singer} where the twistor space $Z$ of a 4-dimensional oriented Riemannian manifold $(M, g)$  is the associated bundle to the principal bundle of orthonormal positive frames with fiber all complex structures in ${\mathbb R}^4$ compatible with the Euclidean scalar product. It has a fiber identified with $\mathbb{S}^2$ and comes equipped with a tautological almost complex structure $J$ which is integrable precisely when $(M,g)$ is self-dual. As Hitchin showed in \cite{hitchin} the only twistor spaces admitting a K\"ahler metric
are $Z \equiv \mathbb{C}P^3$ of $\mathbb{S}^4$ and $Z\equiv Fl_{1,2}$ of $\mathbb{C}P^2$, where $Fl_{1,2}$ is the flag manifold of $\mathbb{C}^3$. This fact makes the general twistor spaces of a self-dual manifolds an appropriate set of examples of non-K\"ahler manifolds. We note here that in this case $Z$ admits a balanced metric \cite{michelson, Oleg}, but only the K\"ahler ones admit an SKT structure \cite{V}. The $\mathbb{S}^2$ fibers of the projection $Z\rightarrow M $ are rational curves and in \cite{ES} another almost complex structure $J_2$ was considered, which has an opposite sign on these fibers and is never integrable. The authors of \cite{ES} showed that a harmonic map from a Riemann surface into $M$ has a lift to a $J_2$-holomorphic map into $Z$.



The Betti numbers of twistor spaces could be computed by Lerray-Hirsch theorem and in \cite{eastwood-singer} the Fr\"olicher spectral sequence $\{E_r\}$ was studied. It  was proved that $\{E_r\}$ always degenerates at most at the second step and a characterization is given for when it degenerates at the first one. There are, in fact, several examples where it fails to degenerate at the first step; in such cases, the Bott-Chern and Aeppli cohomology groups differ from the Dolbeault ones, and the $\partial \overline\partial$-lemma does not hold.

In this work, we study the Bott-Chern and Aeppli cohomologies of twistor spaces and characterize the validity of the $\partial \overline\partial$-lemma in this context, for which no general results were previously available. We compute the Bott-Chern and Aeppli diamonds of a twistor space (Theorem~\ref{thm:bc-cohomology-twistor}), and we show that the twistor space $Z$ of a compact self-dual four-manifold $M$ satisfies the $\partial \overline\partial$-lemma if and only if
\[
h^{1,1}_{BC}(Z) + h^{1,1}_{A}(Z) = 2(b_+(M) + 1)
\quad \text{and} \quad
h^{1,2}_{BC}(Z) = b_1(M)
\]
(Theorem~\ref{theorem9}). As a consequence, we prove that if $M$ is a compact self-dual simply-connected four-manifold such that its twistor space $Z$ satisfies the $\partial \overline\partial$-lemma, then $M$ must be diffeomorphic to either the $4$-sphere or to a connected sum $\#_k \mathbb{C}P^2$ for some $k \geq 1$. Notably, each of these manifolds admits a self-dual metric whose corresponding twistor space satisfies the $\partial \overline\partial$-lemma. Using results from \cite{campana} and \cite{campana2}, we conclude that the twistor space $Z$ of a compact self-dual simply-connected four-manifold $M$ satisfies the $\partial \overline\partial$-lemma if and only if $Z$ lies in Fujiki's class~$\mathcal{C}$, or equivalently, if $Z$ is Moishezon.

In Section~3, we also consider the case of twistor spaces $Z$ of fake projective planes $M$, showing that the Fr\"olicher spectral sequence of $Z$ does not degenerate at the first step; hence, in particular, $Z$ does not satisfy the $\partial \overline\partial$-lemma. Finally, in the last section, we compute explicitly the Dolbeault cohomology of the twistor space $Z$ associated to the flat four-dimensional torus.

\medskip

\noindent{\sl Acknowledgments.} The authors are grateful to the anonymous reviewer for the insightful comments that improved the presentation of the paper.
  This work arose from the AIM workshop \lq \lq Geometric partial differential equations from unified string theories”.  The authors  are grateful for the
great working environment provided by  AIM.  A. Fino, N. Tardini and A. Tomassini are partially supported by the Project PRIN 2022 ``Real and Complex Manifolds: Geometry and Holomorphic Dynamics 2022AP8HZ9'' and by GNSAGA of INdAM. G. Grantcharov is partially supported by a grant from
the Simons Foundation (\#853269). L. Vezzoni is  partially supported by the project PRIN 2022 ``Differential-geometric aspects of manifolds via Global Analysis" and by GNSAGA of INdAM.  We also thank M. Verbitsky for pointing out  to us the reference \cite{campana2} and L.  Di Cerbo for the information on fake projective planes.

\section{Preliminaries}

\subsection{Bott-Chern cohomology and $\partial\overline\partial$-lemma}

Let $X$ be a compact complex manifold of complex dimension $n$ and let $A^{p,q}(X)$ denote the space of complex $(p,q)$-forms on $X$. 
Then the complex \emph{de Rham}, \emph{Dolbeault} and \emph{conjugate Dolbeault}
cohomology groups of $X$ are defined, respectively, as
$$
H^\bullet_{dR}(X;\mathbb{C}):=\frac{\Ker d}{\Imm d}\;,\qquad
H^{\bullet,\bullet}_{\overline\partial}(X):=\frac{\Ker \overline\partial}{\Imm \overline\partial}\;,\qquad
H^{\bullet,\bullet}_{\partial}(X):=\frac{\Ker \partial}{\Imm \partial}\;.
$$
Nevertheless there is no natural map between the de Rham and Dolbeault cohomologies, in this sense a bridge between them is given by
the \emph{Bott-Chern} and the \emph{Aeppli}  cohomology groups that are defined by
$$
H^{\bullet,\bullet}_{BC}(X):=\frac{\Ker \partial\cap\Ker \overline\partial}{\Imm \partial\overline\partial}\;,\qquad
H^{\bullet,\bullet}_{A}(X):=\frac{\Ker \partial\overline\partial}{\Imm \partial+\Imm\overline\partial}.
$$
Notice that
the product induced by the wedge product on forms induces a structure of algebra for the Bott-Chern cohomology 
$H^{\bullet,\bullet}_{BC}(X)$ and a structure of $H^{\bullet,\bullet}_{BC}(X)$-module for the Aeppli cohomology
$H^{\bullet,\bullet}_{A}(X)$.\\
In \cite{schweitzer}, see also \cite{kodaira-spencer}, Hodge theory for the Bott-Chern and the Aeppli cohomologies is developed. In particular, once fixed an Hermitian metric $g$ on $X$ the Bott-Chern and the Aeppli cohomology groups of $X$ are, respectively, isomorphic to the kernel of the following $4^{th}$-order elliptic differential operators
$$ 
\Delta_{BC}^g \;:=\;
\left(\del\delbar\right)\left(\del\delbar\right)^*+\left(\del\delbar\right)^*\left(\del\delbar\right)+
\left(\delbar^*\del\right)\left(\delbar^*\del\right)^*+\left(\delbar^*\del\right)^*\left(\delbar^*\del\right)+\delbar^*\delbar+\del^*\del $$
and
$$\Delta_{A}^g \;:=\; \del\del^*+\delbar\delbar^*+\left(\del\delbar\right)^*\left(\del\delbar\right)+\left(\del\delbar\right)\left(\del\delbar\right)^*+\left(\delbar\del^*\right)^*\left(\delbar\del^*\right)+\left(\delbar\del^*\right)\left(\delbar\del^*\right)^*\,.
$$
Therefore the Bott-Chern and Aeppli cohomologies, such as the de Rham and Dolbeault cohomologies, are finite-dimensional vector spaces. 
We will denote with $b_{\bullet}$ the Betti numbers, $h^{\bullet,\bullet}_{BC}=\dim_{\C} H^{\bullet,\bullet}_{BC}(X)$ and similarly for the other cohomologies.
Moreover, 
the Hermitian duality does not
preserve these cohomologies. Indeed, when an Hermitian metric is fixed on $X$, the Hodge-$*$-operator induces an  isomorphism between the Bott-Chern cohomology and the Aeppli cohomology, 
$$
*:H^{p,q}_{BC}(X)\longrightarrow H^{n-p,n-q}_{A}(X)\,.
$$
As a consequence, we do not have symmetry with respect to the center in the Bott-Chern (and Aeppli) diamond.
More precisely, we have the following equalities:
$h^{p,q}_{BC}(X)=h^{q,p}_{BC}(X)=
h^{n-q,n-p}_{A}(X)=h^{n-p,n-q}_{A}(X)$,
where the first one and the last one follow from the fact that the conjugation preserves
the Bott-Chern and the Aeppli cohomologies respectively (giving a symmetry in the Bott-Chern diamond with respect to the central column).\\

By definition, the identity induces natural maps between the Bott-Chern, Dolbeault, de Rham, and Aeppli cohomologies:
$$\xymatrix{
  & H^{\bullet,\bullet}_{BC}(X) \ar[d]\ar[ld]\ar[rd] & \\
  H^{\bullet,\bullet}_{\partial}(X) \ar[rd] & H^{\bullet}_{dR}(X;\mathbb{C}) \ar[d] & H^{\bullet,\bullet}_{\overline\partial}(X) \ar[ld] \\
  & {\phantom{\;.}} H^{\bullet,\bullet}_{A}(X) \;. &
} $$
Recall that a compact complex manifold is said to satisfy the {\em $\partial\overline\partial$-Lemma} if the natural vertical map $H^{\bullet,\bullet}_{BC}(X)\longrightarrow H^{\bullet,\bullet}_{A}(X)$ is injective. This is equivalent to any of the above maps being an isomorphism (\cite[Lemma 5.15]{deligne-griffiths-morgan-sullivan}).
By \cite{deligne-griffiths-morgan-sullivan} every compact K\"ahler manifold,
and more in general manifolds in class $\mathcal{C}$ of Fujiki (i.e., bimeromorphic to a K\"ahler manifold), satisfy the $\partial\overline\partial$-lemma.
Hence the Bott-Chern and Aeppli cohomologies could provide more informations on a compact complex manifold which does not admit any K\"ahler metric.
Moreover, the $\partial\overline\partial$-lemma can be numerically characterized in terms of dimensions of cohomology groups as follows.

\begin{thm}[\cite{angella-tomassini, angella-tardini}]\label{thm:numbers-deldelbar-lemma}
 Let $X$ be a compact complex manifold. Then,
 $$
 \Delta^k:=\sum_{p+q=k} \left( h^{p,q}_{BC} + h^{p,q}_{A} \right) - 2\, b_k \;\geq\; 0 \;.
 $$
 Moreover, the following facts are equivalent:
 \begin{itemize}
 \item[$1)$] the $\del\delbar$-lemma holds on $X$;
 \medskip
 \item[$2)$] $\Delta^k\;=\; 0 \;$, for any $k\in\mathbb{N}$;
 \medskip
 \item[$3)$] $h^{p,q}_{BC} = h^{p,q}_{A}$, for any $p,q\in\mathbb{N}$.
 \end{itemize}
\end{thm}

\subsection{Twistor spaces}
Let $(M,g)$ be an oriented Riemannian $4$-dimensional manifold. The rank $6$ vector bundle  of $2$-forms $\Lambda^2M$ on $M$ decomposes as the direct sum of two rank 3 vector bundles
$$
\Lambda^2M=\Lambda_+\oplus\Lambda_-
$$
where $\Lambda_\pm$ are the $\pm 1$-eigenspaces of the Hodge-$*$-operator acting on $\Lambda^2M$. The sections of $\Lambda_+$ are called \emph{self-dual $2$-forms} and the sections of $\Lambda_-$ are called \emph{anti-self-dual $2$-forms}.\\
The Riemannian curvature tensor can be thought of as an operator $R : \Lambda^2M \to \Lambda^2M$
called the Riemannian curvature operator. Such operator decomposes under the action of $SO(4)$ as
$$
R=\frac{s}{6}\text{Id}+W_++W_-+\tilde r,
$$
where $W_{\pm}$ are trace-free endomorphisms of $\Lambda_\pm$, and they are called the self-dual
and anti-self-dual components of the Weyl curvature operator. The scalar curvature $s$ acts by scalar multiplication and $\tilde r$ is the trace-free Ricci curvature operator. Recall that an oriented Riemannian $4$-manifold $(M,g)$ is said to be \emph{self-dual} if the Weyl tensor $W=W_+$, i.e., $W_-=0$ and this definition is conformally invariant.\\
The \emph{twistor space} $Z=S(\Lambda_-)$ of a conformal Riemannian manifold
 $\left(M, \left[g\right]\right)$ 
 is the total space of the unit sphere bundle of the vector bundle of anti-self-dual $2$-forms.
This is an $\mathbb{S}^2$-bundle over $M$ whose fiber $C_p$, over a point $p\in M$, is the set of all almost complex structure $J: T_pM\to T_pM$ such that $J$ is orthogonal and compatible with the orientation.
Atiyah, Hitchin, and Singer in \cite{atiyah-hitchin-singer} showed that $Z$ is naturally equipped with
an almost complex structure that is integrable if and only if $W_- = 0$.\\
From now on $M$ will denote a connected smooth compact oriented self-dual $4$-dimensional manifold and $Z=S(\Lambda_-M)$ will denote its twistor space.\\
In \cite[p. 658]{eastwood-singer} Eastwood and Singer, using the Leray-Hirsch theorem, compute the Betti numbers of $Z$:
\begin{equation}\label{betti-numbers}
b_1(Z)=b_1(M)\,,\qquad b_2(Z)=b_2(M)+1\,,\qquad b_3(Z)=2b_1(M)\,
\end{equation}
 and again by \cite{eastwood-singer} the
 Hodge diamond of $Z$ is
\begin{equation}\label{hodge-diamond}
\begin{array}{lclclcl}
& & & h^{0,0}_{\overline\partial}=1 & & &\\
& & h^{1,0} _{\overline\partial}=0& & h^{0,1}_{\overline\partial} =b_1& &\\
& h^{2,0}_{\overline\partial}=0 & & h^{1,1}_{\overline\partial}& & h^{0,2}_{\overline\partial}=b_- &\\
h^{3,0}_{\overline\partial}=0 & & h^{2,1}_{\overline\partial} & & h^{1,2}_{\overline\partial}& & h^{0,3}_{\overline\partial}=0\,.\\
& h^{3,1}_{\overline\partial}=b_- & & h^{2,2}_{\overline\partial} & & h^{1,3}_{\overline\partial}=0 &\\
& & h^{3,2}_{\overline\partial}=b_1 & & h^{2,3}_{\overline\partial}=0 & &\\
& & & h^{3,3}_{\overline\partial}=1 & & &
\end{array}
\end{equation}

Notice that in general $Z$ carries a balanced metric (see \cite{michelson, Oleg}), namely an Hermitian metric $\omega$ such that $d\omega^2=0$ but
is not a K\"ahler manifold. 
In fact Hitchin \cite{hitchin} has shown that only two twistor spaces $Z$ admit a K\"ahler metric, namely 
$\mathbb{C}P^3$ and the space of flags $Fl_{1,2}$ in $\mathbb{C}^3$, that correspond to $\mathbb{S}^4$ and $\mathbb{C}P^2$ respectively. In \cite{campana} Campana has shown that if $Z$ is in class $\mathcal{C}$ of Fujiki (i.e., $Z$ is bimeromorphic to a K\"ahler manifold) then $M$ is either $S^4$ or homeomorphic to a connected sum of copies of $\mathbb{C}P^2$.

In the cases where $Z$ is K\"ahler or in class $\mathcal{C}$ of Fujiki then it satisfies the $\partial\overline\partial$-lemma but, in general it turns out that $Z$ from the cohomological point of view does not have to satisfy the $\partial\overline\partial$-lemma necessarily. Indeed there are examples where the Fr\"olicher spectral sequence does not even degenerate at the first step \cite{eastwood-singer}. For this reason the Bott-Chern and Aeppli cohomology groups are not isomorphic to the Dolbeault cohomology groups.
Notice, however, that in general by  \cite{eastwood-singer} the Fr\"olicher spectral sequence always degenerates at most at the second step.

\section{Bott-Chern and Aeppli cohomologies of twistor spaces and $\partial\overline\partial$-lemma}

In this section we will describe the Bott-Chern and Aeppli cohomologies of the twistor space of a compact self-dual $4$-manifold and we characterize the validity of the $\partial\overline\partial$-lemma. In order to do that we will start with a general lemma.

\begin{lemma}\label{lemma:general}
Let $X$ be a compact complex manifold. Then, the following sequences are exact:
\begin{itemize}
\item[-] $0\to H^{p,0}_{BC}(X)\to H^{p,0}_{\overline\partial}(X)$; \vspace{8pt}
\item[-] $0\to H^{0,q}_{\overline\partial}(X)\to H^{0,q}_{A}(X)\xrightarrow{\overline\partial} H^{0,q+1}_{\partial}(X)$
\end{itemize}
\end{lemma}

\begin{proof}
We first want to show that for any $p$, the map $H^{p,0}_{BC}(X)\to H^{p,0}_{\overline\partial}(X)$ is injective. Notice that $H^{p,0}_{BC}(X)=\text{Ker}\,\partial\cap\text{Ker}\,\overline\partial\cap A^{p,0}(X)$ and $H^{p,0}_{\overline\partial}(X)=\text{Ker}\,\overline\partial\cap A^{p,0}(X)$. The injectivity is then straightforward.
Now we show that the map $\iota:H^{0,q}_{\overline\partial}(X)\to H^{0,q}_{A}(X)$ is injective for any $q$. Indeed, let $[\alpha]_{\overline\partial}\in H^{0,q}_{\overline\partial}(X)$ such that $\iota[\alpha]_{\overline\partial}=[\alpha]_A=0$. Since $H^{0,q}_A(X)=\frac{\text{Ker}\,\partial\overline\partial}{\text{Im}\,\overline\partial}$, then  $\alpha$ is $\overline\partial$-exact and so $[\alpha]_{\overline\partial}=0$. Finally, consider the map $H^{0,q}_{A}(X)\xrightarrow{\overline\partial} H^{0,q+1}_{\partial}(X)=\text{Ker}\,\partial\cap A^{0,q+1}(X)$, defined by $\overline\partial[\alpha]_A=[\overline\partial\alpha]_{\partial}$. Then it is clear that $\text{Im}\,\iota=\text{Ker}\,\overline\partial$.
\end{proof}

Now we can use this lemma to obtain the following

\begin{lemma}\label{lemma:BC}
Let $M$ be a compact self-dual $4$-manifold and let $Z$ be its twistor space. Then, for $p=1,2,3$, the $(p,0)$ and $(0,p)$ Bott-Chern cohomology groups are as follows
$$
H^{p,0}_{BC}(Z)=H^{0,p}_{BC}(Z)=\left\lbrace 0\right\rbrace\,.
$$
\end{lemma}

\begin{proof}
Since by (\ref{hodge-diamond}), $H^{p,0}_{\overline\partial}(Z)=\left\lbrace 0\right\rbrace$ for $p=1,2,3$, it follows by Lemma \ref{lemma:general} that $H^{p,0}_{BC}(Z)\simeq H^{0,p}_{BC}(Z)=\left\lbrace 0\right\rbrace\,$.
\end{proof}

\begin{lemma}\label{lemma:A}
Let $M$ be a compact self-dual $4$-manifold and let $Z$ be its twistor space. Then, for $q=1,2,3$, the $(q,0)$ and $(0,q)$ Aeppli cohomology groups are as follows
$$
H^{q,0}_{A}(Z)=H^{0,q}_{A}(Z)\simeq H^{0,q}_{\overline\partial}(Z)\,.
$$
In particular, $H^{3,0}_{A}(Z)= H^{0,3}_{A}(Z)=\left\lbrace 0\right\rbrace$.
\end{lemma}

\begin{proof}
Notice that by (\ref{hodge-diamond}), for $r=0,1,2,3$, $H^{0,r+1}_{\partial}(Z)=\overline{H^{r+1,0}_{\overline\partial}(Z)}=\left\lbrace 0\right\rbrace$. Hence, applying Lemma \ref{lemma:general} we obtain isomorphisms $H^{q,0}_{A}(Z)=H^{0,q}_{A}(Z)\simeq H^{0,q}_{\overline\partial}(Z)$
\end{proof}

We can now collect all the previous results in the following Theorem.

\begin{thm}\label{thm:bc-cohomology-twistor}
Let $M$ be a compact self-dual $4$-manifold and let $Z$ be its twistor space. Then, the Bott-Chern and Aeppli cohomology groups are as follows
\begin{itemize}
\item[-] $H^{1,0}_{BC}(Z)\simeq H^{0,1}_{BC}(Z)=\left\lbrace 0\right\rbrace$;\vspace{6pt}
\item[-]  $H^{1,0}_{A}(Z)\simeq H^{0,1}_{A}(Z)\simeq H^{0,1}_{\overline\partial}(Z)$;\vspace{6pt}
\item[-]  $H^{2,0}_{BC}(Z)\simeq H^{0,2}_{BC}(Z)=\left\lbrace 0\right\rbrace$;\vspace{6pt}
\item[-]  $H^{2,0}_{A}(Z)\simeq H^{0,2}_{A}(Z)\simeq H^{0,2}_{\overline\partial}(Z)$;\vspace{6pt}
\item[-]  $H^{3,0}_{BC}(Z)\simeq H^{0,3}_{BC}(Z)=\left\lbrace 0\right\rbrace$;\vspace{6pt}
\item[-]  $H^{3,0}_{A}(Z)\simeq H^{0,3}_{A}(Z)=\left\lbrace 0\right\rbrace$.
\end{itemize}
\end{thm}

As a consequence, the Bott-Chern and Aeppli numbers are as follows.

\begin{cor}\label{cor:numbers-twistor}
Let $M$ be a compact self-dual $4$-manifold and let $Z$ be its twistor space. Then, we have the following equalities
\begin{itemize}
\item[-]  $h^{1,0}_{A}(Z)=h^{0,1}_{A}(Z)= b_1(M)$;\vspace{6pt}
\item[-]  $h^{2,0}_{A}(Z)=h^{0,2}_{A}(Z)= b_-(M)$;\vspace{6pt}
\item[-]  $h^{2,3}_{BC}(Z)=h^{3,2}_{BC}(Z)= b_1(M)$;\vspace{6pt}
\item[-]  $h^{1,3}_{BC}(Z)=h^{3,1}_{BC}(Z)= b_-(M)$.
\end{itemize}
\end{cor}
\begin{proof}
The first equality $h^{1,0}_{A}(Z)=h^{0,1}_{A}(Z)= b_1(M)$ follows from Lemma \ref{lemma:A} saying that
$$
H^{1,0}_{A}(Z)\simeq H^{0,1}_{A}(Z)\simeq H^{0,1}_{\overline\partial}(Z)
$$
and the fact $h^{0,1}_{\overline\partial}(Z)=b_1(M)$.\\
The second equality $h^{2,0}_{A}(Z)=h^{0,2}_{A}(Z)= b_-(M)$ again follows from Lemma \ref{lemma:A} saying that
$$
H^{2,0}_{A}(Z)\simeq H^{0,2}_{A}(Z)\simeq H^{0,2}_{\overline\partial}(Z)
$$
and the fact $h^{0,2}_{\overline\partial}(Z)=b_-(M)$.\\
The last two equalities follow from the duality between the Aeppli and Bott-Chern cohomologies.
\end{proof}

As a consequence the Bott-Chern diamond of a twistor space $Z$ is the following
$$
\begin{array}{lclclcl}
& & & h^{0,0}_{BC}=1 & & &\\
& & h^{1,0} _{BC}=0& & h^{0,1}_{BC}=0 & &\\
& h^{2,0}_{BC}=0 & & h^{1,1}_{BC}& & h^{0,2}_{BC}=0 &\\
h^{3,0}_{BC}=0 & & h^{2,1}_{BC} & & h^{1,2}_{BC}& & h^{0,3}_{BC}=0\\
& h^{3,1}_{BC}=b_- & & h^{2,2}_{BC} & & h^{1,3}_{BC}=b_- &\\
& & h^{3,2}_{BC}=b_1 & & h^{2,3}_{BC}=b_1 & &\\
& & & h^{3,3}_{BC}=1 & & &
\end{array}
$$
and dually the Aeppli diamond is
$$
\begin{array}{lclclcl}
& & & h^{0,0}_{A}=1 & & &\\
& & h^{1,0} _{A}=b_1& & h^{0,1}_{A}=b_1 & &\\
& h^{2,0}_{A}=b_- & & h^{1,1}_{A}& & h^{0,2}_{A}=b_- &\\
h^{3,0}_{A}=0 & & h^{2,1}_{A} & & h^{1,2}_{A}& & h^{0,3}_{A}=0\,.\\
& h^{3,1}_{A}=0 & & h^{2,2}_{A} & & h^{1,3}_{A}=0 &\\
& & h^{3,2}_{A}=0 & & h^{2,3}_{A}=0 & &\\
& & & h^{3,3}_{A}=1 & & &
\end{array}
$$

Now, setting as recalled in Theorem \ref{thm:numbers-deldelbar-lemma},
$$
\Delta^k:=\sum_{p,q=0}^k\left(h^{p,q}_{BC}+h^{p,q}_A\right)-2b_k,
$$
by \cite{angella-tomassini} on any compact complex manifold we have that $\Delta^k\geq 0$ for every $k$, and moreover $\Delta^k=0$ for every $k$ if and only if the complex manifold satisfies the $\partial\overline\partial$-lemma.

\begin{cor}\label{cor:delta}
Let $M$ be a compact self-dual $4$-manifold and let $Z$ be its twistor space. Then, we have
\begin{itemize}
\item[-]  $\Delta^1=0$;\vspace{6pt}
\item[-]  $\Delta^2=h^{1,1}_{BC}(Z)+h^{1,1}_A(Z)- 2\left(b_+(M)+1\right)$;\vspace{6pt}
\item[-]  $\Delta^3=4\left(h^{1,2}_{BC}(Z)-b_1(M)\right)$.\vspace{6pt}
\end{itemize}
In particular, the following inequalities hold
\begin{itemize}
\item[-]  $h^{1,1}_{BC}(Z)+h^{1,1}_A(Z)\geq 2(b_+(M)+1)$;\vspace{6pt}
\item[-]  $h^{2,1}_{BC}(Z)=h^{1,2}_{BC}(Z)=h^{2,1}_{A}(Z)=h^{1,2}_{A}(Z)\geq b_1(M)$.
\end{itemize}
\end{cor}

\begin{proof}
First of all on $Z$ we have
$$
\Delta^1=h^{1,0}_{BC}+h^{0,1}_{BC}+h^{1,0}_{A}+h^{0,1}_{A}-2b_1
$$
and by Theorem \ref{thm:bc-cohomology-twistor} and Corollary \ref{cor:numbers-twistor} we have $h^{1,0}_{BC}=h^{0,1}_{BC}=0$ and $h^{1,0}_{A}=h^{0,1}_{A}=b_1$ and so $\Delta^1=0$.\\
We compute now $\Delta^2$ on $Z$. First we have
$$
h^{2,0}_{BC}+h^{1,1}_{BC}+h^{0,2}_{BC}+h^{2,0}_{A}+h^{1,1}_{A}+h^{0,2}_{A}=h^{1,1}_{BC}+2h^{0,2}_{\overline\partial}+h^{1,1}_{A}
$$
now since $h^{0,2}_{\overline\partial}(Z)=b_-(M)$ we have
$$
h^{1,1}_{BC}+2h^{0,2}_{\overline\partial}+h^{1,1}_{A}=h^{1,1}_{BC}+2b_-+h^{1,1}_{A}\,.
$$
Hence,
$$
\Delta^2= h^{1,1}_{BC}+2b_-+h^{1,1}_{A}- 2b_2(Z)=h^{1,1}_{BC}+2b_-+h^{1,1}_{A}- 2(b_2(M)+1)=h^{1,1}_{BC}+2b_-+h^{1,1}_{A}- 2(b_+(M)+b_-(M)+1)
$$
giving
$$
\Delta^2=h^{1,1}_{BC}(Z)+h^{1,1}_{A}(Z)- 2(b_+(M)+1).
$$
Finally, by by Theorem \ref{thm:bc-cohomology-twistor}, Corollary \ref{cor:numbers-twistor}  and by the symmetries of the Bott-Chern and Aeppli cohomology groups on $Z$ we have
$$
h^{3,0}_{BC}+h^{1,2}_{BC}+h^{2,1}_{BC}+h^{0,3}_{BC}+h^{3,0}_{A}+h^{1,2}_{A}+h^{2,1}_{A}+h^{0,3}_{A}=2(h^{3,0}_{BC}+h^{1,2}_{BC}+h^{2,1}_{BC}+h^{0,3}_{BC})=4h^{1,2}_{BC}.
$$
Therefore,
$$
\Delta^3=4h^{1,2}_{BC}(Z)-2b_3(Z),
$$
and since
$$
b_3(Z)=2b_1(M),
$$
we obtain
$$
\Delta^3=4\left(h^{1,2}_{BC}(Z)-b_1(M)\right).
$$
\end{proof}

In a particular, we can characterize the $\partial\overline\partial$-lemma as follows

\begin{thm} \label{theorem9}
Let $M$ be a compact self-dual $4$-manifold and let $Z$ be its twistor space. Then, $Z$ satisfies the $\partial\overline\partial$-lemma if and only if
$$
h^{1,1}_{BC}(Z)+h^{1,1}_A(Z)=2\left(b_+(M)+1\right)
\qquad
\text{and}
\qquad
h^{1,2}_{BC}(Z)=b_1(M).
$$
\end{thm}

Moreover, we can obtain a second numerical characterization of the  $\partial\overline\partial$-lemma. Indeed, recall that by \cite{angella-tardini} (see Theorem \ref{thm:numbers-deldelbar-lemma}) a compact complex manifold satisfies the $\partial\overline\partial$-lemma if and only if $h^{p,q}_{BC}=h^{p,q}_A$ for every $p,q$. Therefore, we have the following

\begin{thm}
Let $M$ be a compact self-dual $4$-manifold and let $Z$ be its twistor space. Then, $Z$ satisfies the $\partial\overline\partial$-lemma if and only if
$$
h^{1,0}_{A}(Z)=0,
\qquad
h^{1,1}_{BC}(Z)=h^{1,1}_A(Z)=b_+(M)+1
\qquad
\text{and}
\qquad
h^{2,0}_{A}(Z)=0.
$$
\end{thm}

In particular, if a twistor space satisfies the $\partial\overline\partial$-lemma one has that
$$
b_1(Z)=h^{1,0}_{BC}+h^{0,1}_{BC}=0
$$
and so the Betti numbers are
$$
b_1(Z)=0\,,\qquad b_2(Z)=b_2(M)+1\,,\qquad b_3(Z)=0\,.
$$
Similarly, by the previous results it follows that

\begin{cor}
Let $M$ be a compact self-dual $4$-manifold and let $Z$ be its twistor space. Then, $Z$ satisfies the $\partial\overline\partial$-lemma if and only if 
its Bott-Chern diamond is the following 
$$
\begin{array}{lclclcl}
& & & h^{0,0}_{BC}=1 & & &\\
& & h^{1,0} _{BC}=0& & h^{0,1}_{BC}=0 & &\\
& h^{2,0}_{BC}=0 & & h^{1,1}_{BC}= b_+(M)+1& & h^{0,2}_{BC}=0 &\\
h^{3,0}_{BC}=0 & & h^{2,1}_{BC}=0 & & h^{1,2}_{BC}=0& & h^{0,3}_{BC}=0\,.\\
& h^{3,1}_{BC}=0 & & h^{2,2}_{BC}=b_+(M)+1 & & h^{1,3}_{BC}=0 &\\
& & h^{3,2}_{BC}=0 & & h^{2,3}_{BC}=0 & &\\
& & & h^{3,3}_{BC}=1 & & &
\end{array}
$$
\end{cor}


Recall that if $M$ is a compact self-dual simply-connected  4-manifold and $Z$ is its twistor space, then by \cite{campana} and  \cite{campana2}  the following are equivalent:
\begin{itemize}
\item[(i)]  $Z$ is in class $\mathcal{C}$ of Fujiki;
\item[(ii)] $Z$ is Moishezon.
\end{itemize}
In particular, $Z$ satisfies the $\partial\overline\partial$-lemma.
However notice that there are small deformations of $Z$ which do not satisfy (i) and (ii), but they do satisfy the $\partial\overline\partial$-lemma. 
 Note that a compact self-dual simply-connected 4-manifold is homeomorphic to $S^4$ or a connected sum of $k$ copies, $k\geq 1$ of complex projective planes by \cite{donaldson}.

%
%
%
%
%
%

\begin{es} 
We now show that the Fr\"olicher spectral sequence of the twistor space $Z$ of fake projective planes $M$ does not degenerate at the first step and so, in particular, it does not satisfy the $\del\delbar$-lemma. Notice that, on the contrary, the twistor space of the complex projective plane is K\"ahler. 
Recall that a {\em  fake projective plane}  is a compact complex surface with the same Betti numbers as the usual
complex projective plane  but  not isomorphic to the complex projective plane. A fake projective plane has ample
canonical divisor, so it is a smooth (and geometrically connected proper)
surface of general type with geometric genus $p_g = 0$ and self-intersection of
canonical class $K^2 = 9$. Note that the fake projective planes satisfies the conditions of the Donaldson's Theorem, except being simply-connected, although with $b_1=0$.  The existence of a fake projective plane was first proved by Mumford  \cite{mumford} and explicit equations of fake projective planes have been determined in \cite{borisov-keum}.   Fake projective planes have Chern numbers $c_1^2 = 3c_2 = 9$ and are complex
2-ball quotients by Aubin  \cite{aubin} and Yau  \cite{yau}.  These  ball quotients are strongly
rigid by Mostow’s rigidity theorem  \cite{mostow}, i.e., are  determined by their  fundamental
group up to holomorphic or anti-holomorphic isomorphisms.  They have been 
 classified as quotients of the two-dimensional complex ball by
explicitly written co-compact torsion-free arithmetic subgroups of $PU(2, 1)$
by Prasad and Yeung    \cite{prasad-yeung} and Cartwright and Steger \cite{cartwright-steger}. The group $PU(2,1)$ is the automorphism group of the 2-ball as well as the isometry group of the metric with constant holomorphic sectional curvature -1. So the metric on the ball induces on the quotient a  metric  with  constant holomorphic sectional curvature which is Einstein and self-dual. It follows that  the fake projective planes admit a self-dual Einstein metric. Moreover, since they have the same Betti numbers of  the complex projective plane, they have $b_1 =0$  and  $b_2^- =0$.\\
We prove now that $E_1\neq E_{\infty}$. By contradiction assume that $E_1= E_{\infty}$, hence, by (\ref{betti-numbers})
$$
0=b_3(Z)=h^{3,0}_{\delbar}(Z)+h^{2,1}_{\delbar}(Z)+h^{1,2}_{\delbar}(Z)+h^{0,3}_{\delbar}(Z)
$$
and so $h^{2,1}_{\delbar}(Z)=0$.\\
By  \cite[Corollary 4.5]{eastwood-singer} since $M$ is Einstein it is regular (in the sense of \cite[page 661]{eastwood-singer}) and so we can apply \cite[Theorem 4.6]{eastwood-singer}. More precisely, since $M$ is regular and $E_1=E_{\infty}$, by \cite[Theorem 4.6]{eastwood-singer} there is a short exact sequence
$$
0\to H^{2}_+(M,\mathbb{C})\to H^{1,1}_{\delbar}(Z)\to \mathbb{C}\to 0\,.
$$
Hence, $h^{1,1}_{\delbar}(Z)=b_++1=2$.\\
However, by \cite[page 658]{eastwood-singer} the kernel of the map $H^{1,1}_{\delbar}(Z)\to H^{2,1}_{\delbar}(Z)$ is canonically identified with $\mathbb{C}\oplus H^{2}_-(M,\mathbb{C})$ and so it has dimension $b_-+1=1$. But this is absurd since, by $h^{2,1}_{\delbar}(Z)=0$ such map is trivial $H^{1,1}_{\delbar}(Z)\to H^{2,1}_{\delbar}(Z)=\left\lbrace 0\right\rbrace$ and the kernel should have dimension $h^{1,1}_{\delbar}(Z)=2$.
\end{es} 

\begin{rmk}
Recall that by \cite[Corollary 5.7]{eastwood-singer}, for any $g \geq 0$, there are metrics on $M = \sharp_g(\mathbb{S}^1 \times\mathbb{S}^3)$ such that the Fr\"olicher spectral sequence of the associated twistor space degenerates at the first step.
We can therefore compute the Hodge numbers of the twistor space using (\ref{hodge-diamond}). Indeed, $h^{0,1}_{\delbar}=b_1=1$,
$$
1=b_2(M)+1=b_2(Z)=h^{2,0}_{\delbar}+h^{1,1}_{\delbar}+h^{0,2}_{\delbar}=h^{1,1}_{\delbar},
$$
and so $h^{1,1}_{\delbar}=1$.
Finally,
$$
2=2b_1(M)=b_3(Z)=h^{3,0}_{\delbar}+h^{2,1}_{\delbar}+h^{1,2}_{\delbar}+h^{0,3}_{\delbar}=2h^{2,1}_{\delbar},
$$
and so $h^{2,1}_{\delbar}=h^{1,2}_{\delbar}=1$.
\end{rmk}

\section{Dolbeault cohomology of the twistor space of the torus}

In this section we consider the twistor space $Z$ of the flat $4$-dimensional torus, which is known to not satisfy the $\partial\overline{\partial}$ lemma, and we compute explicitly the Dolbeault cohomology. The calculations could be generalized to higher-dimensional tori.


 The twistor space $Z$ is the quotient of the bundle
  ${\mathfrak t}^{1,0}\otimes {\it O}(1)$ on $\C P^1\equiv \mathbb{S}^2$  where
  ${\mathfrak t}^{1,0}$ is the $(1,0)$-part of the complexification of
  the $4$-dimensional Abelian algebra $\mathfrak{t}$  (\cite[Ex. 13.64, 13.66]{besse}). More precisely,
 let $[\lambda_1, \lambda_2]$ be the homogeneous coordinates on $\C P^1$.
  On $U_1=\{\lambda\in\C P^1: \lambda_1\neq 0\}$, define $\nu=\frac{\lambda_2}{\lambda_1}$.
  On $U_2=\{\lambda\in\C P^1: \lambda_2\neq 0\}$, define $\mu=\frac{\lambda_1}{\lambda_2}$.

As described in \cite{besse}, on $U_2$,  the space of (0,1)-forms
  \begin{equation}\label{fourteen}
  {\overline\sigma}_1=\frac{\overline{\mu} d\overline{z_1}-d{z}_2}{1+|\mu |^2},
  \hspace{.2in}
  {\overline\sigma}_2=\frac{\overline{\mu} d\overline{z_2}+d{ z}_1}{1+| \mu |^2}, 
  \hspace{.2in}
  {\overline d\mu}
  \end{equation}\label{differential-sigma}
 spans $T^{(0,1)}Z$ at each point. Moreover, the forms ${\overline{\sigma}_1, \overline{\sigma}_2}$
  are holomorphic because
 \begin{equation}
  d{\overline\sigma}_1=\frac{1}{1+|\mu |^2}
      \left( d\overline{\mu}\wedge{\sigma}_2
        -\overline{\mu} d{\mu}\wedge \overline{\sigma}_1
      \right),
  \hspace{.2in}
  d{\overline\sigma}_2=\frac{1}{1+|\mu |^2}
      \left(-\overline{\mu} d{\mu}\wedge{\overline\sigma}_2
       -d\overline{\mu}\wedge{\sigma}_1 \right)
\end{equation}
  are type (1,1)-forms.

As described in \cite{grantcharov-pedersen-poon} on $Z$ one has
  $$H^{k}(Z, p^*{\mathcal{O}} (\ell ))= \mathfrak{t}^{*(0,k)}\otimes S^{\ell+k}\C^{2}$$
 for $k=0,1,2$, where $p$ is the projection from $Z$ to $\mathbb{CP}^1$, and the corresponding Dolbeault groups $H^k(Z, {\mathcal{O}})=H^{0,k}_{\overline\partial}(Z)$ can be explicitly described.
 Define
  \begin{equation}\label{omegas}
  \overline{\Omega}_1=\frac
  {{\overline\lambda}_1d{\overline z}_1-{\overline\lambda}_2dz_2}
  {|\lambda_1|^2+|\lambda_2|^2},
  \hspace{.2in}
  \overline{\Omega}_2=\frac
  {{\overline\lambda}_1d{\overline z}_2+{\overline\lambda}_2dz_1}
  {|\lambda_1|^2+|\lambda_2|^2}.
  \end{equation}
  Then $\{\lambda_1\overline{\Omega}_1, \lambda_2\overline{\Omega}_1,
  \lambda_1\overline{\Omega}_2, \lambda_2\overline{\Omega}_2\}$ form a basis for
   the space $H^1(Z, {\mathcal{O}} )=H^{0,1}_{\overline\partial}(Z)$. 
   Also $H^2(Z, {\mathcal{O}})=H^{0,1}_{\overline\partial}(Z)$, is spanned by
  the twisted 2-forms
  \begin{equation}\label{basic}
    \left(\lambda_1^{2-l}\lambda_2^l\right)
    \overline{\Omega}_1\wedge\overline{\Omega}_2
  \end{equation}
  for $l=0,1,2$

   The space $H^1(Z, {\mathcal{O}})=H^{0,1}_{\overline\partial}(Z)$ has an alternative description.
  For $k=0,1,2,3$ and
  over $p^{-1}({\vec{a}})$, define 1-forms
  \begin{equation}\label{twisted forms}
  {\overline\omega}_k=I_kdx-iI_{\vec{a}}I_kdx,
  \end{equation}
   where $I_0 = Id$, and $I_1,I_2,I_3$ form the standard hypercomplex structure, so $dz_1 = dx+iI_1dx$ and $dz_2 = I_2dx=iI_3dx$, and $I_{\vec{a}} = a_1I_1+a_2I_2+a_3I_3$ for $\vec{a}=(a_1.a_2.a_3)\in S^2$. Note that the norms of $\omega_k$ are all equal to $\sqrt{2}$. Then 
  the space $H^{0,1}_{\overline\partial}(Z)$ is spanned by $\overline{\omega}_k$
  because
  \begin{eqnarray}\label{forms}
  \overline{\omega}_0=\mu {\overline\sigma}_1+{\overline\sigma}_2        =\lambda_1\overline{\Omega}_1+\lambda_2\overline{\Omega}_2, & &
  \overline{\omega}_1=i(\mu {\overline\sigma}_1 -{\overline\sigma}_2)
          =i(\lambda_1\overline{\Omega}_1-\lambda_2\overline{\Omega}_2),
  \nonumber\\
  \overline{\omega}_2=\mu {\overline\sigma}_2-{\overline\sigma}_1
          =\lambda_1\overline{\Omega}_2-\lambda_2\overline{\Omega}_1, & &
  \overline{\omega}_3=i({\overline\sigma}_1+\mu {\overline\sigma}_2)
          =i(\lambda_1\overline{\Omega}_2+\lambda_2\overline{\Omega}_1).
  \end{eqnarray}


The Dolbeault cohomology of $Z$ is described in the following
  
  \begin{thm}\label{thm:dolb-cohom-torus}
 Let $M$ be the $4$-dimensional flat torus and $Z$ its twistor space. Then the following holds.
 \begin{itemize}
 \item[1.] The Hodge diamond of $Z$ is
 $$
\begin{array}{lclclcl}
& & & h^{0,0}_{\overline\partial}=1 & & &\\
& & h^{1,0} _{\overline\partial}=0& & h^{0,1}_{\overline\partial} =4& &\\
& h^{2,0}_{\overline\partial}=0 & & h^{1,1}_{\overline\partial}=4& & h^{0,2}_{\overline\partial}=3&\\
h^{3,0}_{\overline\partial}=0 & & h^{2,1}_{\overline\partial}=4 & & h^{1,2}_{\overline\partial}=4& & h^{0,3}_{\overline\partial}=0\,.\\
& h^{3,1}_{\overline\partial}=3 & & h^{2,2}_{\overline\partial}=4 & & h^{1,3}_{\overline\partial}=0 &\\
& & h^{3,2}_{\overline\partial}=4 & & h^{2,3}_{\overline\partial}=0 & &\\
& & & h^{3,3}_{\overline\partial}=1 & & &
\end{array}
$$
  \item[2.] The non-trivial Dolbeault cohomology groups of $Z$ are
  \begin{itemize}
    \vspace{10pt}
 \item[-] $
  H^{0,1}_{\overline\partial}(Z)=\mathbb{C}\left\langle [\lambda_1\bar\Omega_1]\,,[\lambda_1\bar\Omega_2]\,, [\lambda_2\bar\Omega_1]\,,[\lambda_2\bar\Omega_2]\right\rangle\,;
  $
  \vspace{10pt}
   \item[-] $
  H^{0,2}_{\overline\partial}(Z)=\mathbb{C}\left\langle [\lambda_1\bar\Omega_1\wedge\lambda_1\bar\Omega_2]\,, [\lambda_2\bar\Omega_1\wedge\lambda_2\bar\Omega_2]\,,[[\lambda_1\bar\Omega_1\wedge\lambda_2\bar\Omega_2]\right\rangle\,;
  $
    \vspace{10pt}
 \item[-]  $ 
   \begin{aligned}
  H^{1,1}_{\overline\partial}(Z)= &  \mathbb{C}\left\langle \left[\frac{d\mu\wedge d\bar\mu}{(1+|\mu|^2)^2}\right]\,,
    \left[ (|\lambda_1|^2+|\lambda_2|^2)\,\Omega_1\wedge\bar\Omega_2\right]\,,\right.\\
      &\left. \left[ (|\lambda_1|^2+|\lambda_2|^2)\,\Omega_2\wedge\bar\Omega_1\right]\,,
    \left[ (|\lambda_1|^2+|\lambda_2|^2)\,\left(\Omega_1\wedge\bar\Omega_1-\Omega_2\wedge\bar\Omega_2\right)\right]
  \right\rangle \,;
  \end{aligned}
  $
    \vspace{10pt}
   \item[-] $
      \begin{aligned}
  H^{1,2}_{\overline\partial}(Z)=&\mathbb{C}\left\langle
   [\left(|\lambda_1|^2+ |\lambda_2|^2\right)\lambda_1\,\Omega_1\wedge\bar\Omega_1\wedge\bar\Omega_2]\,,
  [\left(|\lambda_1|^2+ |\lambda_2|^2\right)\lambda_2\,\Omega_1\wedge\bar\Omega_1\wedge\bar\Omega_2]\,,\right.\\
  & \left. [\left(|\lambda_1|^2+ |\lambda_2|^2\right)\lambda_1\,\Omega_2\wedge\bar\Omega_1\wedge\bar\Omega_2]\,,
 [\left(|\lambda_1|^2+ |\lambda_2|^2\right)\lambda_2\,\Omega_2\wedge\bar\Omega_1\wedge\bar\Omega_2]
  \right\rangle\,.
      \end{aligned}
  $
 \end{itemize}
 \end{itemize}
  \end{thm}

\begin{proof}
1. We will start by discussing the Hodge numbers of $Z$. Since $b_1(M)=4$ and $b_-(M)=3$ the only missing numbers in the Hodge diamond (\ref{hodge-diamond}) of the twistor space $Z$ are $h^{1,1}_{\overline\partial}$, $h^{2,2}_{\overline\partial}$, $h^{1,2}_{\overline\partial}$ and $h^{2,1}_{\overline\partial}$. Now, by duality we have $h^{1,1}_{\overline\partial}=h^{2,2}_{\overline\partial}$ and $h^{1,2}_{\overline\partial}=h^{2,1}_{\overline\partial}$, so we are left to prove that $h^{1,1}_{\overline\partial}=4$ and $h^{1,2}_{\overline\partial}=4$.\\
For $h^{1,1}_{\overline\partial}$ we can apply Theorem 5.3 in \cite{eastwood-singer} and since all holomorphic vector fields are parallel we obtain
$$
h^{1,1}_{\overline\partial}=b_0(M)+b_+(M)=1+3=4\,.
$$
Moreover, again by \cite[Theorem 5.3]{eastwood-singer}, the Fr\"olicher spectral sequence degenerates of $Z$ at the first step $E_1=E_{\infty}$ hence
$$
b_3(Z)=h^{3,0}_{\overline\partial}+h^{2,1}_{\overline\partial}+h^{1,2}_{\overline\partial}+h^{0,3}_{\overline\partial}=2\,h^{1,2}_{\overline\partial}
$$
and since $b_3(Z)=2\,b_1(M)=8$ we conclude that $h^{1,2}_{\overline\partial}=4$.\\

\vspace{3pt}
2. Now we show explicitly the representatives of the Dolbeault cohomology classes of $Z$.\\
\begin{itemize}
\item[a)] $
  H^{0,1}_{\overline\partial}(Z)=\mathbb{C}\left\langle [\lambda_1\bar\Omega_1]\,,[\lambda_1\bar\Omega_2]\,, [\lambda_2\bar\Omega_1]\,,[\lambda_2\bar\Omega_2]\right\rangle\,.
  $  \vspace{10pt}\\
It is shown in \cite{grantcharov-pedersen-poon} that $\lambda_1\bar\Omega_1\,,\lambda_1\bar\Omega_2\,, \lambda_2\bar\Omega_1\,,\lambda_2\bar\Omega_2$ represent non-trivial cohomology classes.\\
\item[b)] $
  H^{0,2}_{\overline\partial}(Z)=\mathbb{C}\left\langle [\lambda_1\bar\Omega_1\wedge\lambda_1\bar\Omega_2]\,, [\lambda_2\bar\Omega_1\wedge\lambda_2\bar\Omega_2]\,,[[\lambda_1\bar\Omega_1\wedge\lambda_2\bar\Omega_2]\right\rangle\,.
  $  \vspace{10pt}\\
Since $\lambda_1\bar\Omega_1\,,\lambda_1\bar\Omega_2\,, \lambda_2\bar\Omega_1\,,\lambda_2\bar\Omega_2$ are $\overline\partial$-closed $(0,1)$-forms than their wedge products are still $\overline\partial$-closed. We need to show that the cohomology classes are non-trivial (cf. also \cite{grantcharov-pedersen-poon}).
From now on we will fix on $Z$ the Hermitian metric
$$
\omega=\frac{i}{2}\omega_0\wedge\bar\omega_0+\frac{i}{2}\omega_2\wedge\bar\omega_2+\frac{i}{2}\frac{d\mu\wedge d\bar\mu}{(1+|\mu|^2)^2}\,.
$$
We will show that $\lambda_1\bar\Omega_1\wedge\lambda_1\bar\Omega_2$ is $\overline\partial$-harmonic with respect to this metric and so it determines a non-trivial cohomology class.\\
First of all notice that, on $\lambda_2\neq 0$  (similarly for $\lambda_1\neq 0$),
\begin{equation}\label{lambda-omega}
\lambda_1\bar\Omega_1=\mu\bar\sigma_1\,,\qquad
\lambda_1\bar\Omega_2=\mu\bar\sigma_2\,,\qquad
\lambda_2\bar\Omega_1=\bar\sigma_1\,,\qquad
\lambda_2\bar\Omega_2=\bar\sigma_2\,.
\end{equation}
Applying the $\mathbb{C}$-anti-linear Hodge-$*$-operator we have
$$
*\left(\lambda_1\bar\Omega_1\wedge\lambda_1\bar\Omega_2\right)=f(\mu,\bar\mu)\sigma_1\wedge\sigma_2\wedge\frac{d\mu\wedge d\bar\mu}{(1+|\mu|^2)^2},
$$
where $f(\mu,\bar\mu)$ is a function depending on $\mu,\bar\mu$. It follows by formulas (\ref{differential-sigma}) that
$$
\overline\partial\left(*\left(\lambda_1\bar\Omega_1\wedge\lambda_1\bar\Omega_2\right)\right)=0
$$
and so the $(0,2)$-form is $\overline\partial$-harmonic.\\
For the other classes in $H^{0,2}_{\overline\partial}(Z)$ the argument is similar.\\
 \item[c)]  $   
   \begin{aligned}
  H^{1,1}_{\overline\partial}(Z)= &  \mathbb{C}\left\langle \left[\frac{d\mu\wedge d\bar\mu}{(1+|\mu|^2)^2}\right]\,,
    \left[ (|\lambda_1|^2+|\lambda_2|^2)\,\Omega_1\wedge\bar\Omega_2\right]\,,\right.\\
      &\left. \left[ (|\lambda_1|^2+|\lambda_2|^2)\,\Omega_2\wedge\bar\Omega_1\right]\,,
    \left[ (|\lambda_1|^2+|\lambda_2|^2)\,\left(\Omega_1\wedge\bar\Omega_1-\Omega_2\wedge\bar\Omega_2\right)\right]
  \right\rangle \,.
  \end{aligned} \vspace{10pt}\\
  $
Since the $(1,1)$-form $\frac{d\mu\wedge d\bar\mu}{(1+|\mu|^2)^2}$ is the Fubini-Study metric on $\mathbb{C}P^1$ then it represents a non-trivial cohomology class.\\
Using formulas (\ref{lambda-omega}) and (\ref{differential-sigma}) it is a direct computation to show that the forms
$$
\begin{aligned}
(|\lambda_1|^2+|\lambda_2|^2)\,\Omega_1\wedge\bar\Omega_2&=&
\bar\lambda_1\Omega_1\wedge\lambda_1\bar\Omega_2+\bar\lambda_2\Omega_1\wedge\lambda_2\bar\Omega_2\\
(|\lambda_1|^2+|\lambda_2|^2)\,\Omega_2\wedge\bar\Omega_1&=&
\bar\lambda_1\Omega_2\wedge\lambda_1\bar\Omega_1+\bar\lambda_2\Omega_2\wedge\lambda_2\bar\Omega_1\\
(|\lambda_1|^2+|\lambda_2|^2)\,\left(\Omega_1\wedge\bar\Omega_1-\Omega_2\wedge\bar\Omega_2\right)&=&
\bar\lambda_1\Omega_1\wedge\lambda_1\bar\Omega_1+\bar\lambda_2\Omega_1\wedge\lambda_2\bar\Omega_1\\
&& -\bar\lambda_2\Omega_2\wedge\lambda_2\bar\Omega_2-\bar\lambda_1\Omega_2\wedge\lambda_1\bar\Omega_2
\end{aligned}
$$
are $\overline\partial$-closed.\\
For explanatory purposes, we will compute $\overline\partial\left(\bar\lambda_1\Omega_1\wedge\lambda_1\bar\Omega_2+\bar\lambda_2\Omega_1\wedge\lambda_2\bar\Omega_2\right)$ recalling that $\lambda_1\bar\Omega_2$ and $\lambda_2\bar\Omega_2$ are $\overline\partial$-closed.
Using formulas (\ref{lambda-omega}) we have to compute
$\overline\partial\left(\bar \mu\sigma_1\wedge\mu\bar\sigma_2+\sigma_1\wedge\bar\sigma_2\right)$.
$
\begin{aligned}
\overline\partial\left(\bar \mu\sigma_1\wedge\mu\bar\sigma_2+\sigma_1\wedge\bar\sigma_2\right)=&d\bar\mu\wedge\sigma_1\wedge\mu\bar\sigma_2
+\frac{\bar\mu}{1+\vert\mu\vert^2}\left(d\mu\wedge\bar\sigma_2-\mu d\bar\mu\wedge\sigma_1\right)\wedge\mu\bar\sigma_2+\\
&+\frac{1}{1+\vert\mu\vert^2}\left(d\mu\wedge\bar\sigma_2-\mu d\bar\mu\wedge\sigma_1\right)\wedge\bar\sigma_2=\\
=& \left(\mu-\frac{\vert\mu\vert^2\mu}{1+\vert\mu\vert^2}\right)d\bar\mu\wedge\sigma_1\wedge\bar\sigma_2-\frac{\mu}{1+\vert\mu\vert^2}d\bar\mu\wedge\sigma_1\wedge\bar\sigma_2=0\,.
\end{aligned}
$

\noindent
The differential of the other forms can be computed similarly.

Moreover, with the same argument used above, they are $\overline\partial$-harmonic. Indeed, their Hodge-$*$-operator is of the form
$$
f(\mu,\bar\mu)\alpha\wedge\frac{d\mu\wedge d\bar\mu}{(1+|\mu|^2)^2},
$$
where $f(\mu,\bar\mu)$ is a function depending on $\mu,\bar\mu$ and $\alpha$ is a $(1,1)$-form depending on $\mu,\bar\mu$, $\sigma_i$, $\bar\sigma_j$ and so by formulas (\ref{differential-sigma}) it turns out that $f(\mu,\bar\mu)\alpha\wedge\frac{d\mu\wedge d\bar\mu}{(1+|\mu|^2)^2}$ is $\overline\partial$-closed.\\
 \item[d)] $
      \begin{aligned}
  H^{1,2}_{\overline\partial}(Z)=&\mathbb{C}\left\langle
   [\left(|\lambda_1|^2+ |\lambda_2|^2\right)\lambda_1\,\Omega_1\wedge\bar\Omega_1\wedge\bar\Omega_2]\,,
  [\left(|\lambda_1|^2+ |\lambda_2|^2\right)\lambda_2\,\Omega_1\wedge\bar\Omega_1\wedge\bar\Omega_2]\,,\right.\\
  & \left. [\left(|\lambda_1|^2+ |\lambda_2|^2\right)\lambda_1\,\Omega_2\wedge\bar\Omega_1\wedge\bar\Omega_2]\,,
 [\left(|\lambda_1|^2+ |\lambda_2|^2\right)\lambda_2\,\Omega_2\wedge\bar\Omega_1\wedge\bar\Omega_2]
  \right\rangle\,.
      \end{aligned}
  $ \vspace{10pt}\\
Using again formulas (\ref{lambda-omega}) and (\ref{differential-sigma}) it is a direct computation to show that the forms
$$
\begin{aligned}
\left(|\lambda_1|^2+ |\lambda_2|^2\right)\lambda_1\,\Omega_1\wedge\bar\Omega_1\wedge\bar\Omega_2&=
\bar\lambda_1\Omega_1\wedge\lambda_1\bar\Omega_1\wedge\lambda_1\bar\Omega_2+
\bar\lambda_2\Omega_1\wedge\lambda_1\bar\Omega_1\wedge\lambda_2\bar\Omega_2
\\
\left(|\lambda_1|^2+ |\lambda_2|^2\right)\lambda_2\,\Omega_1\wedge\bar\Omega_1\wedge\bar\Omega_2&=
\bar\lambda_1\Omega_1\wedge\lambda_1\bar\Omega_1\wedge\lambda_2\bar\Omega_2+
\bar\lambda_2\Omega_1\wedge\lambda_2\bar\Omega_1\wedge\lambda_2\bar\Omega_2
\\
\left(|\lambda_1|^2+ |\lambda_2|^2\right)\lambda_1\,\Omega_2\wedge\bar\Omega_1\wedge\bar\Omega_2&=
\bar\lambda_1\Omega_2\wedge\lambda_1\bar\Omega_1\wedge\lambda_1\bar\Omega_2+
\bar\lambda_2\Omega_2\wedge\lambda_1\bar\Omega_1\wedge\lambda_2\bar\Omega_2
\\
\left(|\lambda_1|^2+ |\lambda_2|^2\right)\lambda_2\,\Omega_2\wedge\bar\Omega_1\wedge\bar\Omega_2&=
\bar\lambda_1\Omega_2\wedge\lambda_1\bar\Omega_1\wedge\lambda_2\bar\Omega_2+
\bar\lambda_2\Omega_2\wedge\lambda_2\bar\Omega_1\wedge\lambda_2\bar\Omega_2
\end{aligned}
$$
are $\overline\partial$-closed.\\
For explanatory purposes, we will show that $\bar\lambda_1\Omega_1\wedge\lambda_1\bar\Omega_1\wedge\lambda_1\bar\Omega_2+
\bar\lambda_2\Omega_1\wedge\lambda_1\bar\Omega_1\wedge\lambda_2\bar\Omega_2$ is $\overline\partial$-closed.
Indeed, notice that 
$$
\bar\lambda_1\Omega_1\wedge\lambda_1\bar\Omega_1\wedge\lambda_1\bar\Omega_2+
\bar\lambda_2\Omega_1\wedge\lambda_1\bar\Omega_1\wedge\lambda_2\bar\Omega_2=
-\left(\bar\lambda_1\Omega_1\wedge\lambda_1\bar\Omega_2+\bar\lambda_2\Omega_1\wedge\lambda_2\bar\Omega_2\right)\wedge\lambda_1\bar\Omega_1
$$
and by the computations of $H^{1,1}_{\overline\partial}(Z)$ and $H^{0,1}_{\overline\partial}(Z)$ we have that both $\bar\lambda_1\Omega_1\wedge\lambda_1\bar\Omega_2+\bar\lambda_2\Omega_1\wedge\lambda_2\bar\Omega_2$ and $\lambda_1\bar\Omega_1$ are $\overline\partial$-closed.\\
The differential of the other forms can be computed similarly.

Moreover, with the previous argument , they are $\overline\partial$-harmonic. Indeed, their Hodge-$*$-operator is of the form
$$
f(\mu,\bar\mu)\alpha\wedge\frac{d\mu\wedge d\bar\mu}{(1+|\mu|^2)^2},
$$
where $f(\mu,\bar\mu)$ is a function depending on $\mu,\bar\mu$ and $\alpha$ is a $(1,0)$-form depending on $\mu,\bar\mu$, $\sigma_1$ and $\sigma_2$ and so by formulas (\ref{differential-sigma}) it turns out that $f(\mu,\bar\mu)\alpha\wedge\frac{d\mu\wedge d\bar\mu}{(1+|\mu|^2)^2}$ is $\overline\partial$-closed.\\

\end{itemize}
\end{proof}

\begin{rmk}
Notice that the $(2,2)$-form
$$
\beta=(|\lambda_1|^2+|\lambda_2|^2)^2\,\Omega_1\wedge\bar\Omega_1\wedge\Omega_2\wedge\bar\Omega_2
$$
is $\overline\partial$-closed because it is the product of two $\overline\partial$-closed $(1,1)$-forms. Moreover, it is $\overline\partial$-harmonic because its Hodge-$*$-operator is of the form
$$
f(\mu,\bar\mu)\frac{d\mu\wedge d\bar\mu}{(1+|\mu|^2)^2}
$$
which is still $\overline\partial$-closed. More precisely, $*\beta$ is a $(1,1)$-form and it is the harmonic representative in the cohomology class of the Fubini-Study metric 
$\left[\frac{d\mu\wedge d\bar\mu}{(1+|\mu|^2)^2}\right]$.
\end{rmk}

Concerning the Aeppli and Bott-Chern numbers, the following inequalities hold.
 \begin{prop}
 Let $M$ be the $4$-dimensional flat torus and $Z$ its twistor space. Then the following holds.
 \begin{itemize}
 \item[-] $h^{1,2}_A=h^{2,1}_A=h^{1,2}_{BC}=h^{2,1}_{BC}\geq 4$;
 \vspace{7pt}
 \item[-] $h^{1,1}_A=h^{2,2}_{BC}\geq 5$;
  \vspace{7pt}
 \item[-] $h^{2,2}_A=h^{1,1}_{BC}\geq 4$.
 \end{itemize}
   \vspace{7pt}
 In particular, $\Delta^2>0$.
 \end{prop}
 
 \begin{proof}
 In order to prove that $h^{1,2}_A\geq 4$ notice that the four representatives in Theorem \ref{thm:dolb-cohom-torus} for $H^{1,2}_{\overline\partial}(Z)$ are also non-trivial representatives for classes in $H^{1,2}_{A}(Z)$. Indeed, such forms are $\overline\partial$-closed, and so $\partial\overline\partial$-closed, and moreover their Hodge-$*$-duals (with respect to the metric considered in the proof of Theorem Theorem \ref{thm:dolb-cohom-torus}) is $d$-closed. Hence, such forms are Aeppli-harmonic. In particular, this is showing that the natural map induced by the identity $H^{1,2}_{\overline\partial}(Z)\to H^{1,2}_{A}(Z)$ is injective.\\
With the same argument the natural map induced by the identity $H^{1,1}_{\overline\partial}(Z)\to H^{1,1}_{A}(Z)$ is injective, and so $h^{1,1}_A\geq 4$.
In fact, $h^{1,1}_A\geq 5$, indeed by Theorem \ref{thm:dolb-cohom-torus} we have that
$$
\frac{d\mu\wedge d\bar\mu}{(1+|\mu|^2)^2}\,,\quad
(|\lambda_1|^2+|\lambda_2|^2)\,\Omega_1\wedge\bar\Omega_2\,,\quad
(|\lambda_1|^2+|\lambda_2|^2)\,\Omega_2\wedge\bar\Omega_1
$$
give three representatives for $H^{1,1}_{A}(Z)$, but by direct computations we have that
$$
 \bar\lambda_1\Omega_1\wedge\lambda_1\bar\Omega_1
-\bar\lambda_2\Omega_2\wedge\lambda_2\bar\Omega_2\,,\quad
\bar\lambda_2\Omega_1\wedge\lambda_2\bar\Omega_1
-\bar\lambda_1\Omega_2\wedge\lambda_1\bar\Omega_2
$$
are $\del\delbar$-closed, hence giving two more representatives for $H^{1,1}_{A}(Z)$.\\
 Now we show that $h^{1,1}_{BC}\geq 4$ giving therefore that $h^{2,2}_A\geq 4$. By Theorem \ref{thm:dolb-cohom-torus} we have that
$$
   \begin{aligned}
  H^{1,1}_{\overline\partial}(Z)= &  \mathbb{C}\left\langle \left[\frac{d\mu\wedge d\bar\mu}{(1+|\mu|^2)^2}\right]\,,
    \left[ (|\lambda_1|^2+|\lambda_2|^2)\,\Omega_1\wedge\bar\Omega_2\right]\,,\right.\\
      &\left. \left[ (|\lambda_1|^2+|\lambda_2|^2)\,\Omega_2\wedge\bar\Omega_1\right]\,,
    \left[ (|\lambda_1|^2+|\lambda_2|^2)\,\left(\Omega_1\wedge\bar\Omega_1-\Omega_2\wedge\bar\Omega_2\right)\right]
  \right\rangle \,.
  \end{aligned} \vspace{10pt}\\
  $$
Clearly, the form $\frac{d\mu\wedge d\bar\mu}{(1+|\mu|^2)^2}$ is $d$-closed, but the same holds for the other forms. Indeed, calling $\eta_1$, $\eta_2$, $\eta_3$ the three other representatives respectively 
we have that
$$
\delbar\eta_1=\delbar\eta_2=\delbar\eta_3=0
$$
and
$$
\del\eta_1=-\overline{\delbar\eta_2}=0\,,\quad \del\eta_2=-\overline{\delbar\eta_1}=0\,,\quad \del\eta_3=-\overline{\delbar\eta_3}=0\,.
$$
Hence, they are $d$-closed and $\delbar*$-closed, so $\del\delbar*$-closed, giving therefore non-trivial classes in $H^{1,1}_{BC}(Z)$.\\
Finally, by Corollary \ref{cor:delta} we have
$$
\Delta^2=h^{1,1}_{BC}(Z)+h^{1,1}_A(Z)- 2\left(b_+(M)+1\right)\geq 4+5-8=1>0\,.
$$
 \end{proof}

\end{document}